\documentclass[11pt,reqno]{amsart}

\usepackage{amsmath, amsfonts, amsthm, amssymb,graphicx, color, enumitem}
\newtheorem{theorem}{Theorem}[section]

\theoremstyle{definition}

\theoremstyle{remark}
\newtheorem{remark}{Remark}[section]

\numberwithin{equation}{section}
\allowdisplaybreaks

\renewcommand{\H}{{\bf H}}

\newcommand{\R}{{\mathbb R}}

\def\f{\frac}
\renewcommand{\O}{\Omega}

\def\hf1{^\f{1}{1-\xi^2}}

\def\be{\begin{equation}}
\def\en{\end{equation}}
\def\bs{\begin{split}}
\def\es{\end{split}}
\def\ba{\begin{align}}
\def\ea{\end{align}}

\def\dtre{\partial_3}

\def\dt{\partial_t}

\def\dt{\partial_t}

\def\H{\mathcal H }

\def\hf{\hat F}

\def\cH
{{\mathcal H}}

\author[P. Secchi]{Paolo Secchi}
\address{INdAM Unit \& Department of Civil, Environmental, Architectural Engineering and Mathematics (DICATAM), University of Brescia, Via Valotti 9, 25133 Brescia, Italy}
\email{paolo.secchi@unibs.it}

\title[MHD with perfectly conducting boundary conditions]
{A note on Ideal Magneto-Hydrodynamics with perfectly conducting boundary conditions in the quarter space}


\begin{document}
\dedicatory{\it Dedicated to Prof. Filippo Gazzola on occasion of his 60th birthday}
\bigskip

\begin{abstract}
We consider the initial-boundary value problem in the quarter space for the system of equations
of ideal
Magneto-Hydrodynamics for compressible fluids with perfectly conducting wall boundary conditions. 
On the two parts of the boundary the solution satisfies different boundary conditions, which make the problem an initial-boundary value problem with non-uniformly characteristic boundary.

We identify a subspace ${{\mathcal H}}^3(\Omega)$ of the Sobolev space $H^3(\Omega)$, obtained by addition of suitable boundary conditions on one portion of the boundary, such that for initial data in ${{\mathcal H}}^3(\Omega)$ there exists a solution in the same space ${{\mathcal H}}^3(\Omega)$, for all times in a small time interval. This yields the well-posedness of the problem  combined with a persistence property of full $H^3$-regularity, although in general we expect a loss of normal regularity near the boundary.
Thanks to the special geometry of the quarter space the proof easily follows by the \lq\lq reflection technique\rq\rq.
\end{abstract}

\date{\today}

\maketitle



\section{Introduction}
\label{sect1}

We consider the equations of ideal Magneto-Hydrodynamics (MHD) for the
motion of an electrically conducting compressible  fluid, where "ideal" means that
the effect of viscosity and electrical resistivity is neglected (see \cite{freidberg}):
\begin{equation}\label{mhd}
\begin{cases}
\rho_p(\partial_t+u \cdot\nabla)p +\rho \nabla\cdot u =0,\\
\rho (\partial_t+u \cdot\nabla)u +\nabla p +
H \times(\nabla\times
H )=0,\\
\partial_tH -\nabla\times(u\times H) =0,
\\
(\partial_t+u \cdot\nabla)S=0,
\end{cases}
\end{equation}
in $(0,T)\times\Omega$, where $\Omega$ is a domain in $\R^3$; we denote the boundary of $\Omega$
by $\Gamma$.
In \eqref{mhd} the pressure $p =p (t,x)$, the velocity field $u =u (t,x)=
(u _1,u _2,u _3)$, the magnetic
field $H =H (t,x)=(H _1,H _2,H _3)$ and the entropy $S$ are unknown functions of
time $t$ and space variables $x=(x_1,x_2,x_3)$. The density $\rho $ is
given by the
equation of state $\rho =\rho(p,S )$ where
$\rho>0$ and $\partial \rho/\partial p\equiv\rho_p>0$ for all $p$ and $S$.  We denote 
$\partial_t=\partial/\partial t,\partial_i=\partial/\partial
x_i,\nabla=(\partial_1,\partial_2,\partial_3)$ and use the conventional
notations of vector
analysis. We prescribe the initial conditions
\begin{equation}\label{ic}
(p ,u ,H,S )_{|t=0}=(p^0,u^0,H^0,S^0) \quad\rm{in}\,\,\Omega.
\end{equation}
The system \eqref{mhd} is supplemented with the divergence constraint
\begin{equation}\label{divH}
\nabla\cdot H =0
\end{equation}
on the initial data.

Yanagisawa and Matsumura \cite{MR1092572} have investigated the initial-boundary value problem corresponding to perfectly
conducting wall
boundary conditions. 
To explain the details, let us denote by $\nu$ the unit outward normal to $\Gamma$ and set
\[
\Gamma_0=\{x\in\Gamma: \, (H^0\cdot\nu)(x)=0\}, \quad \Gamma_1=\{x\in\Gamma: \, (H^0\cdot\nu)(x)\not=0\}.
\]
Yanagisawa and Matsumura \cite{MR1092572} prove that for a perfectly conducting wall the 
boundary conditions reduce to
\begin{equation}\label{bc}
\begin{array}{ll}
u  \cdot\nu=0,\quad H \cdot\nu=0 \quad\rm{on}\, (0,T)\times\Gamma_0,
\\
u  =0\quad  \quad\rm{on}\, (0,T)\times\Gamma_1.
\end{array}
\end{equation}
The second condition in the first line of \eqref{bc} is just a restriction on the initial datum $H^0$.
Both boundary conditions in \eqref{bc} are maximal non-negative.

In  \cite{MR1092572} it is considered the case when $\Gamma$ consists only of $\Gamma_{0}$ or $\Gamma_{1}$. In both cases the problems can be reduced to initial boundary value problems for quasi-linear symmetric hyperbolic systems with characteristic boundary of constant multiplicity, see \cite{MR1070840,MR1346224,MR1401431,S96MR1405665}. For the case when $\Gamma$ consists only of  $\Gamma_{0}$ see also \cite{MR1658400,MR1346662,MR1314493,MR1941267,zeroMach}. 
In fact, when $\Gamma$ consists only of  $\Gamma_{0}$ the boundary matrix, that is the coefficient of the normal derivative in the differential operator, has constant rank 2 at $\Gamma_{0}$. Because of a
possible loss of regularity in
the normal direction to the boundary, see \cite{MR2289911,MR336093}, in general the solution of such mixed problems
is not in the usual Sobolev space $H^m(\Omega)$, as for the non-characteristic case,
but in the anisotropic weighted Sobolev space $H^m_\ast(\Omega)$.

On the other hand, when $\Gamma$ consists only of  $\Gamma_{1}$ the boundary matrix has constant rank 6 at $\Gamma_{1}$ (recall that the size of the system is 8). Thus the boundary is again characteristic of constant multiplicity and one could expect the loss of normal regularity. Nevertheless, all the normal derivatives of the vector solution can be estimated by using the nonzero part of the boundary matrix, the special structure of the divergence constraint \eqref{divH} and the fact that the equation for the entropy $S$ is a transport equation. This leads to the proof of the full regularity of the solution in the usual space $H^m(\Omega)$. This is similar to the initial boundary value problem  \eqref{mhd}--\eqref{divH} with boundary conditions
\begin{equation}\label{bc2}
\begin{array}{ll}
u  \cdot\nu=0,\quad H \times\nu=g
\end{array}
\end{equation}
and transversality of the magnetic field at the boundary, see Yanagisawa \cite{MR910816}. The result for the case when $\Gamma$ consists only of  $\Gamma_{1}$ was previously obtained by T. Shirota (not published).

If $\Gamma$ consists of both $\Gamma_{0}$ and $\Gamma_{1}$ the problem is an initial boundary value problem with non-uniformly characteristic boundary, that is characteristic of non-constant multiplicity. If the boundary condition is maximal non-negative, like \eqref{bc}, the existence of weak solutions is classical. However, for non-uniformly characteristic boundary, it is well known that in general weak solutions are not necessarily strong. A sufficient condition for weak=strong is given in Rauch \cite{MR1290491}. A general regularity theory for initial boundary value problems with non-uniformly characteristic boundary, even under Rauch's sufficient conditions, is not yet available. For some results about the regularity of solutions see \cite{MR1376499,MR1759800,MR1666202,MR1769182}.

In the present note we show the local in time well posedness of \eqref{mhd}--\eqref{bc} when the space domain $\O$ is the quarter space. Inspired by \cite{MR1346662}, we identify a subspace $\cH^3(\O)$ of $H^3(\Omega)$, obtained by addition of two suitable boundary conditions on $\Gamma_{0}$, such that for initial data in $\cH^3(\O)$ there exists a solution in the same space $\cH^3(\O)$, for all times in a small time interval. This yields the well-posedness of \eqref{mhd}--\eqref{bc} combined with a persistence property of full $H^3$-regularity, although in general we could only expect a $H^3_\ast$-regularity near $\Gamma_{0}$.
Thanks to the special geometry of the quarter space the proof easily follows by the \lq\lq reflection technique\rq\rq. For other applications of this method in a similar context see \cite{MR1835579,MR2056860}; see also \cite{MR4385755} and references thereinto.


\section{Formulation of the problem, notations and main result}

We denote the quarter space by 
\[\Omega^+=\{x=(x_1,x_2,x_3)\in\R^3\,:\, x_1>0,\, x_3>0\},
\]
and decompose its
boundary as $\Gamma=\Gamma_0\cup\overline\Gamma_1$, where we choose
\begin{equation}\label{def-bdry}
\begin{array}{ll}
\Gamma_0=\{x_1>0, x_3=0\}, \qquad \Gamma_1=\{x_1=0, x_3>0\}.
\end{array}
\end{equation}
The unit outward normal to $\Gamma_0$ is $\nu_0=(0,0,-1)$, and the unit outward normal to $\Gamma_1$ is $\nu_1=(-1,0,0)$. Therefore \eqref{bc} can be rewritten as
\begin{equation}\label{bc1}
\begin{array}{ll}
u_3=0,\quad H_3=0 &\quad\rm{on}\, (0,T)\times\Gamma_0,
\\
u  =0\quad  &\quad\rm{on}\, (0,T)\times\Gamma_1.
\end{array}
\end{equation}
Using \eqref{divH} we rewrite \eqref{mhd} into the following form
\begin{equation}\label{mhd2}
\begin{cases}
\rho_p(\partial_t+u \cdot\nabla)p +\rho \nabla\cdot u =0,\\
\rho (\partial_t+u \cdot\nabla)u +\nabla (p+\frac12|H|^2) -
(H\cdot\nabla)H=0,\\
(\partial_t+u \cdot\nabla)H -(H \cdot\nabla)u  + H \nabla\cdot
u =0,
\\
(\partial_t+u \cdot\nabla)S=0,
\end{cases}
\end{equation}
that can be written in the matrix form as
\begin{equation}\label{quasi}
A_0({{U}})\partial_t{{U}}+\sum_{j=1}^3A_j({{U}})\partial_j{{U}}=0
\end{equation}
for $U=(p,u,H,S)^T$, with

\begin{equation}\label{}
\begin{array}{ll}
A_0(U)= \begin{pmatrix}
{\rho_p}/{\rho}&\underline 0&\underline 0&0\\ 
\underline 0^T&\rho I_3 &O_3&\underline 0^T\\ \underline 0^T&O_3&I_3&\underline 0^T\\
0&\underline 0&\underline 0&1
 \end{pmatrix},
\end{array}
\end{equation}
\begin{equation}\label{}
\begin{array}{ll}
A_j(U)=\begin{pmatrix}(\rho_p/\rho) u_j&\delta_j&\underline 0&0
\\
\delta^T_j&\rho
u_jI_3&\delta_j\otimes H-H_jI_3&\underline 0^T\\
\underline 0^T&(\delta_j\otimes H)^T-H_jI_3&u_jI_3&\underline 0^T
\\
0&\underline 0&\underline 0&u_j
\end{pmatrix},
\end{array}
\end{equation}
where
$\delta_j=(\delta_{1j},\delta_{2j},\delta_{3j}),\delta_{kj}$ is the Kronecker delta,
$\delta_j\otimes H$ is the
$3\times 3$ matrix $(\delta_{jk}H_i),k\downarrow1,2,3,i\rightarrow1,2,3$, $\underline 0=(0,0,0)$.
The quasilinear system \eqref{quasi} is symmetric hyperbolic if the state equation $\rho=\rho(p,S)$ satisfies
the hyperbolicity condition $A_0>0$, i.e.
\begin{equation}\label{hyper}
\rho(p,S)>0,\quad \rho_p(p,S)>0.
\end{equation}
%
We define
\begin{equation}\label{defN}
\begin{array}{ll}
N=\{U=(p,u,H,S)^T:\overline{\O^+}\to\R^8\,:\, u_3=H_3=0\quad\text{on }\Gamma_0\},\\
\\
N^\perp=\{U=(p,u,H,S)^T:\overline{\O^+}\to\R^8\,:\, p=u_1=u_2=H_1=H_2=S=0\quad\text{on }\Gamma_0\},
\end{array}
\end{equation}
and introduce the Sobolev subspace
\begin{equation}\label{defH}
\begin{array}{ll}
\H^3(\O^+)=\{U\in H^3(\O^+)\,:\, U\in N,\,\dtre U\in N^\perp, \partial_{33}^2U\in N\}.
\end{array}
\end{equation}

Given the system \eqref{quasi} for $U$ with initial condition $U_{|t=0}=U^0=(p^0,u^0,H^0,S^0)^T$, we recursively define
$\dt^kU^0=(\dt^kp^0,\dt^ku^0,\dt^kH^0,\dt^kS^0)^T$, $k\geq 1$, by formally
taking $k-1$ time derivatives of the equations, solving for
$\partial_t^kU$ and
evaluating it at time $t=0$ in terms of $U^0$ and its space derivatives; for $k=0$ we set $\dt^0U^0=U^0$.

The main result of the paper is given by the following theorem. 
%
The result can be extended to any order $m$ of regularity, by showing the existence of solutions and persistence property of $H^m$-regularity in a suitable subspace $\H^m$, defined in a similar way as $\H^3$ in \eqref{defH}, with the addition of more \lq geometric' properties on derivatives up to order $m-1$.
\begin{theorem}\label{th-main}
Let $\rho\in C^{4}$ and
$U^0=(p^0,u^0,H^0,S^0)^T$ be such that $U^0-(0,\underline 0,c,\underline 0)^T \in  \H^3(\O^+)$ for some constant $c\not=0$, $\rho (p^0,S^0)>0,\rho_p(p^0,S^0)>0$ in
$\overline{\Omega^+}$, $\nabla\cdot H^0=0$ in $\Omega^+$, $H^0_1\not=0$ on $\overline\Gamma_1$. We also assume that the initial datum satisfies the compatibility conditions
\begin{equation}\label{comp-u}
\begin{array}{ll}
\dt^ku^0=0 \quad \text{for} \quad k=0,1,2, \quad \text{on}\,\, \Gamma_1.
\end{array}
\end{equation}
Then there exists $T>0$ such that the mixed problem \eqref{ic}, \eqref{bc1}, \eqref{quasi}
has a unique solution 
$$
U-(0,\underline 0,c,\underline 0)^T\in \cap_{k=0}^3 C^k([0,T];\H^{3-k}_{}(\O^+))
$$
satisfying \eqref{divH}, \eqref{hyper} in $[0,T]\times\overline{\O^+}$.
\end{theorem}
\begin{remark}\label{compa}
The compatibility conditions associated with the boundary conditions on $\Gamma_0$ are 
$$
\dt^ku^0_3=0 \quad\text{for}\quad k=0,1,2,\quad H^0_3=0 \quad\text{on}\quad\Gamma_0.
$$
These compatibility conditions are not explicitly prescribed in the statement of Theorem \ref{th-main}, because they are automatically satisfied if $U^0-(0,\underline 0,c,\underline 0)^T \in  \H^3(\O^+)$. In fact, if $k=0$ we have by definition $U^0-(0,\underline 0,c,\underline 0)^T \in N$, that is $u^0_3=H^0_3=0$ on $\Gamma_0$. If $k=1$ we write 
\[
\dt U^0=-\sum_{j=1}^3\hat{A}_j({{U^0}})\partial_j{{U^0}},
\]
where we have denoted $\hat{A}_j(U')=A_0(U')^{-1}A_j(U')$.
Then $\dt U^0\in N$ on $\Gamma_0$ easily follows from the \lq geometric'  properties
\begin{equation}\label{geometric1}
\begin{array}{ll}
\hat{A}_j(U')N\subset N, \quad \hat{A}_j(U')N^\perp\subset N^\perp, \quad j=1,2,\\
\hat{A}_3(U')N\subset N^\perp, \quad \hat{A}_3(U')N\subset N,
\end{array}
\end{equation}
for all $U'\in N$, see \cite[Section 3]{MR1346662}.
With similar arguments we show that $\dt^2 U^0\in N$ on $\Gamma_0$ using \eqref{geometric1} and
\begin{equation}\label{geometric2}
\begin{array}{ll}
\dt\hat{A}_j(U')N\subset N, \quad \dt\hat{A}_j(U')N^\perp\subset N^\perp, \quad j=1,2,\\
\dt\hat{A}_3(U')N\subset N^\perp, \quad \dt\hat{A}_3(U')N\subset N,
\end{array}
\end{equation}
for all $U'\in N$, see again \cite[Section 3]{MR1346662}.
\end{remark}
\begin{remark}\label{}
From the definition of $\Gamma_1$ and the assumption $H^0_1-c \in  H^3(\O^+)$ in Theorem \ref{th-main} it follows that $H^0_1$ is uniformly bounded away from zero on the unbounded set $\Gamma_1$. This is crucial for the uniform invertibility of the non-singular part of the boundary matrix on $\Gamma_1$ and the proof of full regularity in \cite{MR1092572}, Theorem 2.7. As far as we know, the opposite case, when $c=0$, is open.
\end{remark}

\section{Proof of Theorem \ref{th-main}}

Let us introduce the half-space  
\[\Omega=\{x=(x_1,x_2,x_3)\in\R^3\,:\,  x_1>0\},
\]
whose boundary is 
\[\partial\Omega=\{x=(x_1,x_2,x_3)\in\R^3\,:\,  x_1=0\}.
\]
Given the initial datum $U^0:\O^+\to \R^8$ as in the statement of Theorem \ref{th-main}, we consider the extension 
\begin{equation}\label{extension}
\begin{array}{ll}

\tilde{U}^0=(\tilde{p}^0,\tilde u^0,\tilde H^0,\tilde S^0)^T:\O\to\R^8,
\end{array}
\end{equation}
where $\tilde u_3^0,\tilde H^0_3$ are respectively the odd extension of $ u^0_3, H^0_3$, with respect to $x_3$, and $\tilde{p}^0,\tilde u^0_1$, $\tilde u^0_2$, $\tilde H^0_1,\tilde H^0_2,\tilde S^0$ are respectively the even extension of ${p}^0, u^0_1, u^0_2, H^0_1, H^0_2, S^0$, with respect to $x_3$. For instance,
\begin{equation*}
\begin{array}{ll}\label{}
\tilde u_3^0(x_1,x_2, x_3)=\begin{cases}
u^0_3(x_1,x_2, x_3) &\quad\text{for  }x_3\geq0,\\
-u^0_3(x_1,x_2, -x_3) &\quad\text{for  }x_3<0,
\end{cases}
\end{array}
\end{equation*}
and similarly for $\tilde H^0_3$;
\begin{equation*}
\begin{array}{ll}\label{}
\tilde p^0(x_1,x_2, x_3)=\begin{cases}
p(x_1,x_2, x_3) &\quad\text{for  }x_3\ge0,\\
p(x_1,x_2, -x_3) &\quad\text{for  }x_3<0,
\end{cases}
\end{array}
\end{equation*}
and similarly for $\tilde u^0_1$, $\tilde u^0_2$, $\tilde H^0_1,\tilde H^0_2,\tilde S^0$.
Next, given $\tilde U^0$, we consider the following initial-boundary value problem on $\O$:
\begin{equation}\label{mixedO}
\begin{cases}
A_0({{U}})\partial_t{{U}}+\sum_{j=1}^3A_j({{U}})\partial_j{{U}}=0 \qquad&\text{on}\,\,(0,T)\times\O\,,\\
u=0  \qquad&\text{in}\,\,(0,T)\times\partial\O\,,\\
U_{|t=0}=\tilde U^0 \qquad&\text{in}\,\,\O\,.
\end{cases}
\end{equation}
The existence of the solution to \eqref{mixedO} follows from \cite[Theorem 2.7]{MR1092572}, that we recall here for the reader's convenience, with some small change to adapt it to our notation. We notice that the original version also considers the case of the unbounded domain with compact smooth boundary that we don't need.
\begin{theorem}[\cite{MR1092572}, Theorem 2.7]\label{YM}
Let $\O'$ be an unbounded domain in $\R^3$ with sufficiently smooth and compact boundary $\partial\O'$ with outward unit vector $n$ (respectively a half space $\R^3_+$). Let $m\ge3$ be an integer. Suppose that $U_0'-(c',\underline0,\underline0,0)^T\in H^m(\O')$ for some constant $c'>0$ (respectively $U_0'-(c',\underline0,c,\underline0)^T\in H^m(\R^3_+)$ for some constants $c'>0, c\not=0$) and that $U_0'=(p_0',u_0',H_0',S_0')$ satisfies the conditions
\begin{equation}\label{comp-H}
\begin{array}{ll}
\nabla\cdot H_0'=0, \quad p_0'>0\text{  in  }\O',\quad H_0'\cdot n\not=0 \text{  on  }\partial\O',
\end{array}
\end{equation}
and the compatibility conditions
\begin{equation}\label{comp-cond}
\begin{array}{ll}
\dt^ku_0'=0 \quad\text{for}\quad k=0,\dots,m-1, \text{  on  }\partial\O'.

\end{array}
\end{equation}
Then there exists a constant $T>0$ such that the problem \eqref{mixedO} with initial datum $U_0'$ has a unique solution
$$
U-(c',\underline 0,\underline 0,0)^T\in \cap_{k=0}^m C^k([0,T];H^{m-k}_{}(\O'))
$$
$\Big(respectively \,\,\,  
U-(c',\underline 0,c,\underline 0)^T\in \cap_{k=0}^m C^k([0,T];H^{m-k}_{}(\R^3_+))
\Big)$.
\end{theorem}
In \eqref{comp-cond} the terms $\dt^ku_0'$ are computed in terms of $U_0'$ as explained before for  $\dt^kU^0$ computed in terms of $U^0$.
\begin{remark}\label{pressure}
In \cite{MR1092572} the authors assume for $\rho$ a constitutive law of the form $\rho =\rho(p,S )$ where
$\rho>0$ and $\partial \rho/\partial p\equiv\rho_p>0$ for all $p>0$ and $S$. Accordingly, in  \cite[Theorem 2.7]{MR1092572} they require the initial pressure to be strictly positive $p_0'>0$ in $\O$. This is different from the present paper where we assume $\rho>0$ and $\partial \rho/\partial p\equiv\rho_p>0$ for all $p$ and $S$. It is easily checked that the result of Theorem \ref{YM} holds in our case as well, that is with $c'=0$.
\end{remark}
\begin{remark}[\cite{MR1092572}]\label{}
The assumptions that $\nabla\cdot H_0'=0$ in $\O$ and $H_0'\cdot n\not=0$ on $\partial\O$ in \eqref{comp-H} imply that the boundary $\partial\O$ consists of more than two connected components except when $\O$ is a half space.
\end{remark}

We wish to apply Theorem \ref{YM} for $\O'=\R^3_+$, $m=3$, the initial datum $U_0'=\tilde U^0$ and $c'=0$, as observed in Remark \ref{pressure}. 

Let $U^0$ be as in Theorem \ref{th-main}. It follows from $U^0-(0,\underline 0,c,\underline 0)^T \in  \H^3(\O^+)$ that the extension \eqref{extension} satisfies $\tilde U^0-(0,\underline 0,c,\underline 0)^T \in  H^3(\O)$. Moreover, from $\nabla\cdot H^0=0$ in $\Omega^+$ we readily get $\nabla\cdot \tilde H^0=0$ in $\Omega$. Since we are taking the even extension of $H^0_1$ with respect to $x_3$, from $H^0_1\not=0$ on $\overline\Gamma_1$ it also follows that $\tilde H^0_1\not=0$ on $\partial\O$.

To apply Theorem \ref{YM} it remains to check the compatibility conditions \eqref{comp-cond}. 
If $k=0$, then $u^0=0$ on $\Gamma_1$ gives immediately by $x_3$-reflection that $\tilde u^0=0$ on $\partial\O$. For $k=1$ we observe that by definition
the assigned value of $\dt u^0$ in $\O^+$ is equivalent to saying that $U^0$ formally solves the equation for the velocity $\eqref{mhd2}_2$ in $\O^+$. On the other hand, by direct computation if $U^0$ solves  $\eqref{mhd2}$ in $\O^+$, then $\tilde U^0$ solves  $\eqref{mhd2}$ in $\O$. Then \eqref{comp-u} for $k=1$ gives $\dt\tilde u^0=0$ on $\partial\O$.
A similar argument, which also involves $\dt p^0$, $\dt H^0$, $\dt S^0$, gives $\dt^2\tilde u^0=0$ on $\partial\O$.

We apply Theorem \ref{YM} and obtain the unique solution $U$ to \eqref{mixedO} such that
$$ 
U-(0,\underline 0,c,\underline 0)^T\in \cap_{k=0}^3 C^k([0,T];H^{3-k}_{}(\R^3_+))\,.
$$ 
Now we define
$
\tilde{U}=(\tilde{p},\tilde u,\tilde H,\tilde S)^T
$
where $\tilde u_3,\tilde H_3$ are respectively the odd extension of $ u_3, H_3$ (restricted to $\O^+$) with respect to $x_3$, and $\tilde{p},\tilde u_1$, $\tilde u_2$, $\tilde H_1,\tilde H_2,\tilde S$ are respectively the even extension of ${p}, u_1, u_2, H_1, H_2, S$ (restricted to $\O^+$) with respect to $x_3$.  For instance,
\begin{equation*}
\begin{array}{ll}\label{}
\tilde u_3(t,x_1,x_2, x_3)=\begin{cases}
u_3(t,x_1,x_2, x_3) &\quad\text{for  }x_3\geq0,\\
-u_3(t,x_1,x_2, -x_3) &\quad\text{for  }x_3<0,
\end{cases}
\end{array}
\end{equation*}
and similar definition for $\tilde H_3$;
\begin{equation*}
\begin{array}{ll}\label{}
\tilde p(t,x_1,x_2, x_3)=\begin{cases}
p(t,x_1,x_2, x_3) &\quad\text{for  }x_3\ge0,\\
p(t,x_1,x_2, -x_3) &\quad\text{for  }x_3<0,
\end{cases}
\end{array}
\end{equation*}
and similar definitions for $\tilde u_1$, $\tilde u_2$, $\tilde H_1,\tilde H_2,\tilde S$.

By direct calculations, we prove that $\tilde{U}$
is also a solution to the initial-boundary value problem \eqref{mixedO}. Thus the uniqueness of the
solution of \eqref{mixedO} implies that $U=\tilde{U}$.
This yields that $u_{3}$, $H_{3}$ 
are odd functions in $x_3$, and hence they satisfy the conditions
\begin{equation*}\label{bc3}
\begin{array}{ll}
u_3=0,\quad H_3=0 &\quad\rm{on}\, (0,T)\times\Gamma_0.
\end{array}
\end{equation*}
Therefore $U$ restricted to $(0,T)\times\O^+$ is the desired solution to our initial-boundary value problem
\eqref{ic}, \eqref{bc1}, \eqref{quasi}.
The proof of Theorem \ref{th-main} is complete.

\section{Acknowledgement}
The research was supported in part by the Italian research
project PRIN 20204NT8W4-002.
The author thanks Prof. A. Matsumura for his kind introduction to the problem with characteristic boundary of non constant multiplicity.

\bibliographystyle{plain}

\bibliography{perf-conduct-quarter-space}

\end{document}